# Integrated Energy Exchange Scheduling for Multi-microgrid System with Electric Vehicles

Dai Wang, *Student Member, IEEE*, Xiaohong Guan, *Fellow, IEEE*, Jiang Wu, *Member, IEEE*, Pan Li, Peng Zan and Hui Xu

***Abstract*—Electric vehicles (EVs) can be considered as flexible mobile battery storages in microgrids. For multiple microgrids in an area, coordinated scheduling on charging and discharging are required to avoid power exchange spikes between the multi-microgrid system and the main grid. In this paper, a two-stage integrated energy exchange scheduling strategy for multi-microgrid system is presented, which considers EVs as storage devices. Then several dual variables, which are representative of the marginal cost of proper constraints, are utilized to form an updated price, thereby being a modification on the original electricity price. With this updated price signal, a price-based decentralized scheduling strategy is presented for the Microgrid Central Controller (MGCC). Simulation results show that the two-stage scheduling strategy reduces the electricity cost and avoids frequent transitions between battery charging/discharging states. With the proposed decentralized scheduling strategy, each microgrid only needs to solve its local problem and limits the total power exchange within the safe range.**

## I. NOMENCLATURE

*Index:*

| | |
|---|---|
| $i$ | Index for microgrid, $i = 1,2,\ldots,I$. |
| $j$ | Index for EV, $j = 1,2,\ldots,J_i$. |
| $m$ | Index for EV's trip, $m = 1,2,\ldots,M$. |
| $t$ | Index for hour, $t = 1,2,\ldots,T$. |
| $s$ | Index for scenario, $s = 1,2,\ldots,S$. |

*Parameters:*

| | |
|---|---|
| $\underline{B}_{i,s}(t)$ | Lower bound of the aggregate batteries' remaining energy in microgrid $i$ during $t$ in scenario $s$ (kWh). |
| $\overline{B}_{i,s}(t)$ | Upper bound of the aggregate batteries' remaining energy in microgrid $i$ during $t$ in scenario $s$ (kWh). |
| $C(t)$ | Original electricity price during $t$ (CNY/kWh). |
| $c^{\mathrm{b}}$ | Capital cost of EV battery (CNY). |
| $c^{\mathrm{l}}$ | Cycle life of EV battery (cycle). |
| $d_{i,j,s}^{m}$ | Driving distance of EV $j$ in microgrid $i$ for trip $m$ in scenario $s$ (km). |
| $E$ | Capacity of EV battery (kWh). |
| $\underline{soc}$ | Lower bound of battery state-of-charge (SOC). |
| $k$ | Electric drive efficiency (km/kWh). |
| $M_{i,j,s}$ | Trip number of EV $j$ in microgrid $i$ in scenario $s$. |
| $n_{i,s}(t)$ | Number of EVs connected to the grid in microgrid $i$ during $t$ in scenario $s$. |
| $\overline{p}_{\mathrm{ev}}$ | Maximal charging/discharging power of individual EV (kW). |
| $\overline{P}_{i,s}(t)$ | Maximal charging power of the aggregate battery storage in microgrid $i$ during $t$ in scenario $s$ (kW). |
| $\underline{P}_{i,s}(t)$ | Minimal discharging power of the aggregate battery storage in microgrid $i$ during $t$ in scenario $s$ (kW). |
| $P_{i,t}^{\mathrm{L}}$ | Base load of microgrid $i$ during $t$ (kW). |
| $P_{i,s}^{\mathrm{W}}(t)$ | Wind power in microgrid $i$ during $t$ in scenario $s$ (kW). |
| $t_{i,j,s}^{\mathrm{a}}$ | Arrival time of EV $j$ in microgrid $i$ in scenario $s$. |
| $t_{i,j,s}^{\mathrm{d}}$ | Departure time of EV $j$ in microgrid $i$ in scenario $s$. |
| $\kappa$ | Upper bound of daily cycle number. |
| $\tau_s$ | Probability of scenario $s$. |
| $\rho$ | Scaling factor， representing the cost of increasing the power exchange capacity. |
| $\Delta_{i,s}^{\mathrm{B}}(t)$ | Energy change of the aggregate battery storage caused by EV departure and arrival in microgrid $i$ during $t$ in scenario $s$ (kWh). |
| $\boldsymbol{\Gamma}_{i,s}(t)$ | Arrival vector, $\boldsymbol{\Gamma}_{i,s}(t) = (\Gamma_{i,1,s}(t),\ldots,\Gamma_{i,J,s}(t))$, $\Gamma_{i,j,s}(t) = 1$(or 0) if EV $j$ arrives home (or not) at $t$ in scenario $s$. |
| $\mathbf{soc}_{i,s}^{\mathrm{ini}}$ | Initial SOC vector, $\mathbf{soc}_{i,s}^{\mathrm{ini}} = (soc_{i,1,s}^{\mathrm{ini}}, soc_{i,2,s}^{\mathrm{ini}},\ldots,soc_{i,J,s}^{\mathrm{ini}})$, $soc_{i,j,s}^{\mathrm{ini}}$ is the initial SOC when EV $j$ arrives home in scenario $s$. |
| $\boldsymbol{\Phi}_{i,s}(t)$ | Departure vector, $\boldsymbol{\Phi}_{i,s}(t) = (\Phi_{i,1,s}(t),\ldots,\Phi_{i,J,s}(t))$, $\Phi_{i,j,s}(t) = 1$(or 0) if EV $j$ leaves home (or not) at $t$ in scenario $s$. |
| $\mathbf{soc}_{i,s}^{0}$ | Departure SOC vector, $\mathbf{soc}_{i,s}^{0} = (soc_{i,1,s}^{0},\ldots,soc_{i,J,s}^{0})$, $soc_{i,j,s}^{0}$ is the SOC of EV $j$ when it leaves home in scenario $s$. |

*Variables:*

| | |
|---|---|
| $B_{i,s}(t)$ | Remaining energy of the aggregate battery storage |

The research presented in this paper is supported in part by the National Natural Science Foundation (61221063, U1301254, 61304212), 863 High Tech Development Plan (2012AA011003) and 111 International Collaboration Program, of China.

D. Wang, X. Guan, J. Wu, P. Li and P. Zan are with the Ministry of Education Key Lab for Intelligent Networks and Network Security (MOE KLINNS), School of Electronic and Information Engineering, Xi'an Jiaotong University, Xi'an, China. J. Wu is the corresponding author of this paper. (Email: jwu@sei.xjtu.edu.cn).

H. Xu is with Binhai Power Supply Bureau, Tianjin, China.

in microgrid $i$ during $t$ in scenario $s$ (kWh)

$\overline{C}^*(t)$ Updated price signal during $t$ (CNY/kWh).

$soc_{i,j,s}(t)$ SOC of EV $j$ in microgrid $i$ during $t$ in scenario $s$ (kWh).

$P_{i,s}^{\mathrm{B}}(t)$ Total battery charging (>0) or discharging (<0) power in microgrid $i$ during $t$ in scenario $s$ (kW).

$P^{\mathrm{cap}}$ Upper bound of power exchange between the multi-microgrid system and the distribution grid (kW).

$P_{i,s}^{\mathrm{ex}}(t)$ Power exchange between the distribution grid and the microgrid $i$ during $t$ in scenario $s$ (kW).

$p_{i,j,s}^{\mathrm{c}}(t)$ Charging power of EV $j$ in microgrid $i$ during $t$ in scenario $s$ (kW).

$p_{i,j,s}^{\mathrm{d}}(t)$ Discharging power of EV $j$ in microgrid $i$ during $t$ in scenario $s$ (kW).

$u_{i,j,s}(t)$ Binary variable, "1" if EV $j$ in microgrid $i$ in scenario $s$ is charging during $t$, "0" otherwise.

$v_{i,j,s}(t)$ Binary variable, "-1" if EV $j$ in microgrid $i$ in scenario $s$ is discharging during $t$, "0" otherwise.

## II. Introduction

MICROGRID is a localized entity consisting of electricity generation, energy storage, and loads and connected to a traditional centralized grid [1]. Renewable energy such as photovoltaic generation and wind power may be utilized in a microgrid as distributed energy sources for avoiding the long distance transmission and reducing carbon emission [2]. Storage devices in a micorgrid not only help make better use of renewable energies but also regulate the load peaks [3].

Electric vehicles (EVs) can be viewed as "mobile batteries". Many countries are making great efforts to increase the penetration of EVs into the market [4]. The Chinese government has also laid down policies such as "The 12th Five-Year Plan on Developing EV Technology" to promote the EV industrialization and improve the integration of EVs into power systems. Since household EVs park in work places or at home 22 hours per day on average [5], the microgrid controller can take advantage of idle EVs as battery storages. The investigations reported in [6-8] showed that EVs have strong ability to provide spinning reserve and frequency regulation through vehicle-to-grid (V2G).

Microgrid central controller (MGCC) is responsible for optimizing one or multiple microgrids' operation [9]. All EVs connected to the power grid can be considered as an aggregate battery storage. However, due to dynamic behaviors of EV charging and discharging, the capacity of this aggregate storage is not constant. Driving patterns (departure time, arrival time, number of daily trips and distance of each trip) have great influences on the battery status, and should be taken into account when MGCC decides the energy scheduling.

Many strategies are proposed to reduce the total cost or minimize the power losses within a microgrid [10-17]. Most of them focused on individual microgrid optimization without considering the possible impact on the main grid. For multiple microgrids in a same area, if each microgrid aims to maximize its own profits and manages the energy exchange with the main grid independently under the external given Time-of-Use (TOU) price, all microgrids will have homogenous behaviors. They will charge their storage devices during the price-valley time and discharge their storage devices during the price-peak time. As a result, the total power exchange may reach a new high peak and exceed the allowable limits. Clearly with its role of an aggregator that acts in the interest of multiple microgrids [9], MGCC should consider the impact of a multi-microgrid system on the main grid.

In this paper, we try to address two issues: 1) how to make use of EVs' storage ability for microgrid scheduling with EV charging/discharging behaviors; 2) how to limit the power exchange peak between the multi-microgrid system and the main grid. Following this line, we present a two-stage integrated scheduling strategy for energy exchange to and from a multi-microgrid system. Considering the stochastic characteristics of EVs and volatility of wind power, we apply the scenario tree method to solve the stochastic optimization problem. At the first stage, the MGCC aims to minimize the total electricity cost of the multi-microgrid system and limit the total power exchange within a safe range. At the second stage, the aggregate charging/discharging power is allocated to each EV and charging/discharging transitions are optimized to increase battery life. A decentralized scheduling strategy is derived based on the dual decision variables at the first stage. The key contributions of this paper are summarized as follows: due to the dynamic behaviors of EV charging/discharging, EV driving patterns that have direct influences on the performance of the aggregate battery storage are introduced into the hierarchical scheduling system and the optimization model. A price-based decentralized scheduling strategy is presented based on the dynamic price update. The multi-microgrid system can respond to the price signal and reduce the power exchange peak. The performances are evaluated and compared through a set of numerical simulations.

The rest of this paper is organized as follows: Related publications are reviewed in Section III. The two-stage problem is formulated in Section IV. Then an updated price signal is developed to establish a decentralized energy exchange strategy in Section V. Simulation and analysis are shown in Section VI. Finally, conclusions and future work are drawn in Section VII.

## III. Related Work

The EV charging/discharging allocation is studied in many recent works. In [7], Sortomme et al. presented algorithms to find optimal charging rates with the objective of maximizing the aggregator's profit. The same authors developed an effective V2G algorithm to optimize energy and ancillary services scheduling in [8]. This algorithm maximizes profits to the aggregator while providing additional system flexibility and peak load shaving to the utility and low cost of EV charging to the customer. In [18], the authors explored the relationship between feeder losses, load factor and load variance and proposed three optimal charging control algorithms to minimize the impacts of PEVs charging on the distribution system. In a recent work [19], the authors proposed an online coordination method to optimally charge PEVs in order to maximize the PEV owner's satisfaction and to minimize system operating costs. The proposed charging architecture guarantees the feasibility of the charging decisions by means of a prediction unit that can forecast future EVs power demand and through a two-stage optimization unit. These works provided comprehensive results on the EV charging/discharging and ancillary services, but those authors did not consider the EV's energy storage function in microgrids. Besides, EV driving patterns which directly determine the parameters of the aggregate battery storage should obtain more attention.

Several studies have looked at EVs' storage ability in microgrid or distribution network. Working as "mobile batteries", EVs' charging/discharging behaviors can be controlled to match the volatile wind power within the microgrid whilst meeting EV users' requirement [20, 21]. In [20], the authors exploited three coordinated wind-PEV energy dispatching approaches in the vehicle-to-grid (V2G) context. In [21], Jianhui el al. used a new unit commitment model which can simulate the interactions among plug-in hybrid electric vehicles (PHEVs), wind power, and demand response (DR). The benefits of EVs which work as storage devices are also discussed in different kinds of microgrids, including residential microgrid [22], industrial microgrid [23] and commercial building microgrid [24]. However, the variability of the aggregate battery storage caused by EVs' driving patterns is not discussed. The battery cycle life reduction caused by V2G implementation is not considered either.

There are also some related literatures about energy management strategies designed for MGCC. Among these strategies, genetic algorithm [10] and lagrangian relaxation [11] are proposed to solve the problem of dispatchable generator's output scheduling and the battery charging planning with distributed renewable energy sources from a perspective of centralized control. To avoid the problem that the autonomous behavior of each consumer and the interaction between the consumers and the generators are always ignored in centralized control strategies, game theories [12-14] and pricing mechanisms [15, 16] are exploited to facilitate the microgrid operation. However, all these methods focus on the cost minimization or benefit maximization within a microgrid. They do not consider the possible impact of multiple microgrids in an area when energy storage devices are integrated into the system.

The differences between our work and these related works are summarized as follows: 1) related works [7, 8, 18, 19] on EV charging/discharging allocation do not consider the energy storage function from a MGCC perspective. EV driving patterns which directly determine the parameters of the aggregate battery storage should gain more attention; 2) the variability of aggregate battery caused by individual EV's driving pattern and the battery cycle life reduction caused by V2G implementation are not considered in [20-24]. There are also few works on a proper framework to integrate the coordinated EV charging/discharging allocation into the microgrid energy management; 3) how to restrain the impact of multiple neighboring microgrids on the main grid is another focus of this work, which, however, is not considered by most of related works [10-16].

## IV. Problem Formulation

The multi-microgrid system is an entity, which is formed of multiple individual microgrids. These individual microgrids are geographically close and connected to the same distribution bus through a common transformer. All these microgrids are at the control of a single MGCC. The structure of the multi-microgrid system is shown in Fig.1. The top level is MGCC. It uses the market prices of electricity and grid security concerns to determine the amount of power that each microgrid should draw from the distribution system. According to the definition in [9], the MGCC represents the function of an aggregator or energy service provider that can acts in the interests of multiple microgrids in an area. The second level is local microsource controllers (MC) and load controllers (LC), which receives the control signal from the MGCC and react to that signal correspondingly. The bottom level is EV aggregator which communicates with each EV owner and schedules the charging and discharging. It also uploads the individual information on all EVs, namely departure time, arrival time and driving distance [20, 25, 26].

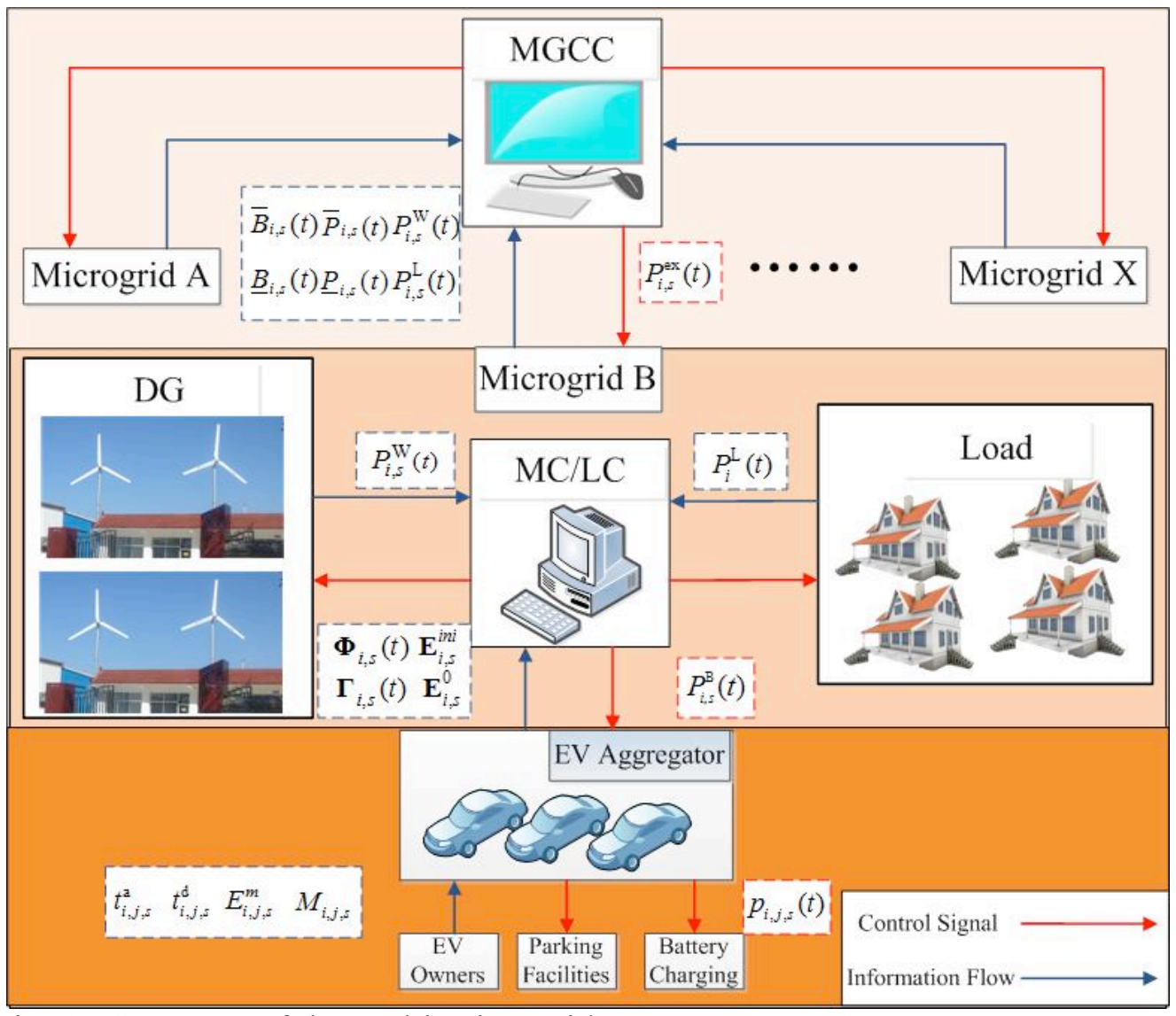

Fig. 1. Structure of the multi-microgrid system.

In this section, a coordinated two-stage energy exchange strategy is proposed. The first stage is MGCC scheduling stage. The MGCC considers all EVs connected to the grid as the aggregate battery storage. It aims at achieving the minimum cost by scheduling the power exchange between the multi-microgrid system and the main distribution grid. The second stage is EV aggregator scheduling stage. Because the first stage only decides the optimal total charging/discharging power, the second stage needs to allocate the total charging/discharging power to each individual EV. In addition, it is the aggregator's duty to minimize the transitions of battery charging/discharging states, since the cycle life of EV battery is limited. It should be noted that we focus on the energy dispatch problem, of which resolution is typically between 15 minutes to 1 hour. We assume that the system dynamics caused by any fluctuation at a finer time scale (several seconds to 5 minutes) will be handled by primary and secondary frequency controls [27].

### A. Stochastic Characteristics of EVs and Wind Power

The charging and discharging behaviors of EVs make them considered as either power loads or power supplies when connected to the power grid. The arrival time is determined from the last trip to home, and the departure time is determined by the first outgoing trip from home next day. Between the arrival time and departure time, EVs are considered at home and charging/discharging is available. With the recognition that the driving patterns of the EV users are with a stochastic nature, we derive the statistical model of EV driving patterns (the arrival time, departure time, number of daily trips and distance of each trip) by using the real-world driving data which is obtained from the NHTS2009 [28].

The departure time distribution is fitted in the form of chi-square distribution.

$$P_{\mathrm{DEP}}(t_{\mathrm{depn},i}) = \frac{t_{\mathrm{depn},i}^{(\upsilon-2)/2} e^{-t_{\mathrm{depn},i}/2}}{2^{\upsilon/2}\Gamma(\upsilon/2)} \tag{1}$$

where $\Gamma(\cdot)$ is defined as $\Gamma(z)=\int_0^\infty t^{z-1}\mathrm{e}^{-1}\mathrm{d}t$ , $t_{\mathrm{depn},i}$ is the normalized departure time at the $i$th departure time window and defined as $t_{\mathrm{dep},i} / \Delta t$ , $t_{\mathrm{dep},i}$ is the departure time at the $i$th departure time window, and $\Delta t$ is the discretized window size. $\upsilon$ is determined to minimize the root-mean-square error of the response variable by applying sequential quadratic programming.

The arrival time at a concerning departure time is expressed as a conditional probability. The distribution of arrival time at the $i$th departure time is:

$$P_{\mathrm{ARR,DEP}}(t_{\mathrm{arr}} \mid t_{\mathrm{dep},i}) = \frac{1}{\sqrt{2\pi\sigma_i^2}} e^{-\frac{(t_{\mathrm{arr}}-\mu_i)^2}{2\sigma_i^2}} \tag{2}$$

Where $t_{\mathrm{arr}}$ is the arrival time, $\mu_i$ is the mean of the arrival time at the $i$th departure time window, and $\sigma_i$ is the standard deviation of the arrival time at the $i$th departure time window.

The distance of each trip conforms to the truncated power-law (3) with exponent $\beta = 1.25$ , $d_0 = 1.8$ and $\alpha = 20$ .

$$P(d) = (d_0 + d)^{-\beta} \exp(-d / \alpha) \tag{3}$$

The distribution of number of daily trips can also easily be obtained from the dataset. Details about the driving pattern models can be seen in our previous work [29]. These statistical models will be used to generate EV driving data. We assume that the EV aggregator can collect or forecast the individual information about all EVs as in [20, 26].

It is known that the forecasted wind power generation $P^{\mathrm{W}}$ is determined by forecasted wind speed $v_{\mathrm{f}}$ and turbine parameters, i.e., cut-in wind speed $v_{\mathrm{ci}}$ , cut-out wind speed $v_{\mathrm{co}}$ , nominal wind speed $v_r$ and nominal power of wind power generator $P^{\mathrm{r}}$ [30, 31]. In this work, we adopt the existing model:

$$P^{\mathrm{W}} = \begin{cases} 0, & v_{\mathrm{f}} < v_{\mathrm{ci}} \text{ or } v_{\mathrm{f}} > v_{\mathrm{co}} \\ \dfrac{v_{\mathrm{f}}^3 - v_{\mathrm{ci}}^3}{v_{\mathrm{r}}^3 - v_{\mathrm{ci}}^3} P^{\mathrm{r}}, & v_{\mathrm{ci}} \le v_{\mathrm{f}} \le v_{\mathrm{r}} \\ P^{\mathrm{r}}, & v_{\mathrm{r}} \le v_{\mathrm{f}} \le v_{\mathrm{co}} \end{cases} \tag{4}$$

Due to the volatility of wind power, it cannot be accurately forecasted. The wind power generation can be assumed to be subject to a normal distribution $N(\mu,\sigma^2)$ with forecasted wind power as its expected value ( $\mu$ ) and a percentage of $\mu$ as its volatility ( $\sigma$ ) [32]. Scenarios are generated according to these distribution models. To alleviate the computational burden caused by a large number of scenarios, the fast forward algorithm [33] is applied to complete the scenario reduction.

## B. MGCC Scheduling Problem

At the MGCC scheduling stage, all EVs connected to the grid can be considered as an aggregate battery storage. MGCC receives the information on the aggregate battery storage, including the bound of battery stored energy ( $\underline{B}_{i,s}(t)$ , $\overline{B}_{i,s}(t)$ ) and the bound of total charging/discharging power ( $\underline{P}_{i,s}(t)$ , $\overline{P}_{i,s}(t)$ ), rather than the detailed information of every EV. In the multi-microgrid system, wind power consumption is assumed at zero marginal cost and is used in priority. Since the wind power cannot cover total energy demand, the remaining load is supplied by the main distribution grid. Under the TOU price, every microgrid tries to minimize its cost, thereby charging its aggregate battery during price-valley time and discharging the battery during the price-peak time. The MGCC pays for the energy consumption from the distribution grid and receives revenues for the energy that is fed back to the distribution grid. Since the batteries are already bought by EV owners, the MGCC can only consider the energy cost. The objective (5) is to minimize the average operation cost of the multi-microgrid system, and obtain the optimal electricity exchange capacity. It is assumed that benefits are shared by all consumers. The problem (Pb-1) is formulated with the scenario tree method, as follows:

$$\min_{P_{i,s}^{\mathrm{ex}}(t), P^{\mathrm{cap}}} \sum_{s=1}^{S} \tau_s \sum_{i=1}^{I} \sum_{t=1}^{T} C_t \cdot P_{i,s}^{\mathrm{ex}}(t)\Delta t + \rho \cdot P^{\mathrm{cap}} \tag{5}$$

$$\text{s.t. } P_{i,s}^{\mathrm{W}}(t) + P_{i,s}^{\mathrm{ex}}(t) = P_i^{\mathrm{L}}(t) + P_{i,s}^{\mathrm{B}}(t) \tag{6}$$

$$B_{i,s}(t) = B_{i,s}(t-1) + P_{i,s}^{B}(t)\Delta t + \Delta_{i,s}^{B}(t) \tag{7}$$

$$\underline{B}_{i,s}(t) \le B_{i,s}(t) \le \overline{B}_{i,s}(t) \tag{8}$$

$$\underline{P}_{i,s}(t) \le P_{i,s,}^{\mathrm{B}}(t) \le \overline{P}_{i,s}(t) \tag{9}$$

$$-\overline{P}_i^{\mathrm{ex}} \le P_{i,s}^{\mathrm{ex}}(t) \le \overline{P}_i^{\mathrm{ex}} \tag{10}$$

$$-P^{\mathrm{cap}} \le \sum_{i=1}^{I} P_{i,s}^{\mathrm{ex}}(t) \le P^{\mathrm{cap}} \tag{11}$$

In this model, $\Delta t$ is the step size of scheduling time. $P^{\mathrm{cap}}$ is a decision variable in the objective function, which is related to the construction cost of infrastructure and facilities. $\rho$ is a scaling factor, which converts the fixed investment to the daily cost during the time horizon of the formulation. The value of $\rho$ can be chosen to adjust the weight of the electricity exchange capacity in the objective against the electricity cost. A relatively large value of $\rho$ may result in a lower electricity exchange capacity, but a higher electricity cost, and vice versa.

Equation (6) denotes the demand and supply balance inside each microgrid. In this paper, it is believed that the load can be accurately forecasted. Thus, $P_{i,s}^{\mathrm{ex}}(t)$ is determined by $P_{i,s}^{\mathrm{B}}(t)$ and vice versa. To clearly illustrate the energy exchange between the multi-microgrid system and the main distribution grid, $P_{i,s}^{\mathrm{ex}}(t)$ is considered as a decision variable. Because the capacity of the aggregate battery fluctuates with the number of EVs connected to the grid, the SOC cannot represent the exact amount energy stored in the battery storage. Instead, we use (7)-(9) to impose practical constraints on the battery storage at the MGCC scheduling stage. Equation (7) updates the amount of energy in battery storage. Equation (8) and (9) ensure that the battery energy and the charging/discharging power lie in feasible ranges, respectively. It should be pointed out that all microgrids are connected to a common low-voltage bus and all energy is exchanged through the common bus. Equation (10) limits the exported or imported power, due to the transmission line capacity between each individual microgrid and the common bus. The security constraint of power exchange between the multi-microgrid system and the main grid is shown in (11). It requires that the total energy exchange should be less than $P^{\text{cap}}$ for all scenarios, even for the worst case.

Since the number of EVs connected to the grid varies during each period, the parameters of the aggregate battery storage are not constant. The acquisition of aggregate battery storage parameters in scenario $s$ are shown in (12)-(17), as follows:

$$\Delta_{i,s}^{\text{B}}(t) = [\mathbf{\Gamma}_{i,s}(t)\mathbf{soc}_{i,s}^{\text{ini}} - \mathbf{\Phi}_{i,s}(t)\mathbf{soc}_{i,s}^{0}]\cdot E \quad (12)$$

$$\underline{B}_{i,s}(t) = [\mathbf{\Phi}_{i,s}(t+1)\mathbf{soc}_{i,s}^{0} + n_{i,s}(t)\underline{soc}\ ]\cdot E \quad (13)$$

$$\overline{B}_{i,s}(t) = n_{i,s}(t)E \quad (14)$$

$$n_{i,s}(t) = n_{i,s}(t-1) + \left\|\mathbf{\Gamma}_{i,s}(t)\right\|_1 - \left\|\mathbf{\Phi}_{i,s}(t)\right\|_1 \quad (15)$$

$$\overline{P}_{i,s}(t) = n_{i,s}(t)\overline{p}_{\text{ev}} \quad (16)$$

$$\underline{P}_{i,s}(t) = -n_{i,s}(t)\overline{p}_{\text{ev}} \quad (17)$$

The energy change of the aggregate battery storage caused by EV arrival and departure is shown in (12). Equation (13) defines the lower bound of remaining battery capacity, which must meet the driving demand and the battery security constraint. The lower bound changes with the number of EVs connected to the grid. It must satisfy the driving demand for the EVs which are leaving next moment. At the same time, the lower bound must also guarantee that the SOC of all connected EVs cannot be lower than the safe limit. Therefore, the lower bound in (13) includes two parts. The first part is the required energy for the EVs which are leaving next moment. The second part is minimum energy to guarantee that the remaining energy of all connected EVs will not be lower than the safe limit $\underline{soc}$. The upper bound of stored energy is shown in (14), which is determined by the number of EVs connected to the grid. The number of EVs at home is calculated by (15). The upper and lower bounds of charging/discharging power are obtained by (16) and (17), respectively.

*Remark*: We provide a way to extend this model to AC network. AC power flow can be calculated by the embedded solver in EPRI's Open-source Distribution System Simulator (OpenDSS) [34]. Without loss of generality, a mathematical formulation (18) can be used to implicitly represent a variety of power flow constraints [35].

$$\mathbf{G}_s(t) \le 0 \quad \forall s, \forall t \quad (18)$$

The bus voltage is limited by (19).

$$\underline{V}_i \le V_{i,s}(t) \le \overline{V}_i \quad \forall s, \forall t \quad (19)$$

where $\underline{V}_i$ and $\overline{V}_i$ are the lower and upper voltage limits for the individual microgrid $i$. The main focus of this work, however, is to exploit the potential of EVs to work as an aggregate battery storage and alleviate the impact of energy exchange peak on the main grid. The goal of this work is to provide insights and lay theoretical foundation for the energy scheduling of multi-microgrid system with EVs. Therefore, we consider the problem from the energy aspect. The following work is discussed using the DC model to avoid distracting readers' attention.

## C. EV Aggregator Scheduling Problem

Once the optimal total charging/discharging power $P_{i,s}^{\text{B}}(t)$ is obtained through the MGCC scheduling problem, the EV aggregator needs to allocate the total power to each individual EV. Therefore, four factors of driving patterns, namely departure time, arrival time, the number of daily trips and the distance of each trip, should be considered. Among these factors, departure time and arrival time determine the feasible charging/discharging period, and the other two factors decide the energy demand.

Technical model:

$$soc_{i,j,s}(t+1) = soc_{i,j,s}(t) + \frac{[p^{\text{c}}{}_{i,j,s}(t) + p^{\text{d}}{}_{i,j,s}(t)]\Delta t}{E} \quad t \in [t_{i,j,s}^{\text{a}}, t_{i,j,s}^{\text{d}}] \quad (20)$$

$$soc_{i,j,s}(t_{i,j,s}^{d}) \ge \frac{\sum_{m=1}^{M} d_{i,j,s}^{m}}{kE} \quad (21)$$

$$v_{i,j,s}(t)\overline{p}_{ev} < p^{\text{d}}{}_{i,j,s}(t) < 0 \quad (22)$$

$$0 < p^{\text{c}}{}_{i,j,s}(t) < u_{i,j,s}(t)\overline{p}_{ev} \quad (23)$$

$$\underline{soc} < soc_{i,j,s}(t) < 100\% \quad (24)$$

$$\sum_{j=1}^{J_i}(p^{\text{c}}{}_{i,j,s}(t) + p^{\text{d}}{}_{i,j,s}(t)) = P_{i,s}^{\text{B}}(t) \quad (25)$$

$$\frac{1}{2}\sum_{t=1}^{T}(u^{\text{z}}{}_{i,j,s}(t) - v^{\text{z}}{}_{i,j,s}(t)) \le \kappa \quad (26)$$

$$u_{i,j,s}(t) - v_{i,j,s}(t) \le 1 \quad (27)$$

The technical part of EV aggregator scheduling problem consists of the dynamic state of charging equation and a set of battery parameter constraints, displayed in (20)-(27). SOC dynamics are shown in (20). Equation (19) defines the minimum SOC when the EV leaves home. The limit of charging/discharging power is shown in (22) and (23). The depth of discharge (DOD) must be less than ($1-\underline{soc}$) to avoid over-discharge, as shown in (24). Equation (25) means that the sum of charging/discharging power of all EVs should be equal to the optimal aggregate value that is achieved from the MGCC scheduling problem. The daily number of cycles is limited by (26). Equation (27) ensures that the EV cannot charge and discharge simultaneously.

Economic model:

The aggregator's total cost calculations $K^{\text{total}}$ are divided into energy cost $K_{i,s}^{\text{e}}$ and battery degradation cost $K_{i,s}^{\text{w}}$.

$$K_{i,s}^{\text{total}} = K_{i,s}^{\text{e}} + K_{i,s}^{\text{w}} \quad (28)$$

The electricity cost $K_{i,s}^{\text{e}}$ is determined by the optimal total charging/discharging power $P_{i,s}^{\text{B}}(t)$, which is obtained through the MGCC scheduling problem.

$$K_{i,s}^{\text{e}} = \sum_{t=1}^{T} P_{i,s}^{\text{B}}(t)\Delta t \cdot C(t) \quad (29)$$

Practically, the battery degradation is related to various factors, including ambient temperature, depth of discharge, charging/discharging rate and etc [36, 37]. It is difficult to get an accurate analytical equation to define the influence of partial charging/discharging on battery's cycle life [38]. Their large uncertainties and non-linear correlations would not yield better solutions and their consideration would increase very much the computation time and the complexity of the problem. In this paper, we define the combination of one charging and one discharging as a cycle and roughly believe the life of EV battery just depends on the number of charging/discharging cycles. Therefore, the battery degradation cost is indicated by charging/discharging cycles [3].

$$K_{i,s}^{\mathrm{w}} = \sum_{t=1}^{T}\sum_{j=1}^{J_i}\frac{1}{2}(u^{\mathrm{z}}{}_{i,j,s}(t) - v^{\mathrm{z}}{}_{i,j,s}(t)) \cdot \frac{c^b}{c^l} \tag{30}$$

From another perspective, since the $u^{\mathrm{z}}{}_{i,j,s}(t)$ and $v^{\mathrm{z}}{}_{i,j,s}(t)$ are variables in equation (30), the objective of aggregator can be interpreted as minimizing state transitions between charging and discharging, which is obviously harmful to the EV battery. Moreover, it is assumed that at least one charging/discharging cycle should be required due to the daily driving demand. If the aggregator can schedule the coordinated charging/discharging processes within this inherent cycle, it is believed that no extra degradation is introduced by the coordinated scheduling. Assistant equations (31)-(32) help to calculate the state transitions between charging and discharging.

$$\begin{cases} u_{i,j,s}^{\mathrm{z}}(t) - v_{i,j,s}^{\mathrm{z}}(t) \le 1 \\ u_{i,j,s}^{\mathrm{z}}(t) - v_{i,j,s}^{\mathrm{za}}(t) \le 1 \\ u_{i,j,s}^{\mathrm{za}}(t) - v_{i,j,s}^{\mathrm{z}}(t) \le 1 \end{cases} \tag{31}$$

$$\begin{cases} u_{i,j,s}(t) - u_{i,j-1,s}(t) = u_{i,j,s}^{\mathrm{z}}(t) + v_{i,j,s}^{\mathrm{za}}(t) \\ v_{i,j,s}(t) - v_{i,j-1,s}(t) = v_{i,j,s}^{\mathrm{z}}(t) + u_{i,j,s}^{\mathrm{za}}(t) \end{cases} \tag{32}$$

where $u_{i,j,s}^{\mathrm{za}}(t) \in \{0,1\}$ and $v_{i,j,s}^{\mathrm{za}}(t) \in \{-1,0\}$ are the assistant parameters to balance (32) when the state of the battery changes from charging (discharging) to idle, since this state transition does not generate penalty.

In sum, the aggregator's scheduling problem is:

$$\min_{p^{\mathrm{c}}{}_{i,j,s}(t), p^{\mathrm{d}}{}_{i,j,s}(t), u_{i,j,s}(t), v_{i,j,s}(t)} K_{i,s}^{\mathrm{total}} \tag{33}$$

s.t. (20)-(27), (31)-(32)

## V. Price-based Decentralized Scheduling Strategy

Since multiple microgrids are coupled by (11), the scheduling in the MGCC stage is a centralized scheme which ensures that the total power exchange will not exceed the limits. Under the centralized scheme, the MGCC minimize the global cost and have to directly control all devices. Nevertheless, such centralized scheduling strategy does not consider the autonomous behaviors of each microgrid, thus cannot guarantee the economic fairness for each microgrid. In this section, we are trying to develop a decentralized scheduling strategy based on a dynamic price update. This price signal uses the optimal values of dual variables, which are representative of the marginal cost of proper constraints. It is calculated one time (day ahead) and will be broadcasted to each individual microgrid as an internal price. Since this price is a modification on the original TOU price, we call it "updated price signal". Each microgrid can adjust its own energy management based on this updated price signal. It ensures that the total power exchange will stay within the safe range.

From MGCC scheduling formulation, we see that $P_{i,s}^{\mathrm{ex}}(t)$, $P_{i,s}^{\mathrm{B}}(t)$ and $B_{i,s}(t)$ are variables. We can standardize the linear problem as follows:

$$\min_{\mathbf{X}} \mathbf{C} \cdot \mathbf{X} \tag{34}$$

$$\text{s.t. } \mathbf{A} \cdot \mathbf{X} \ge \mathbf{D} \tag{35}$$

where $\mathbf{X} = \left[\mathbf{B};\ \mathbf{P}^{\mathrm{B}};\ \mathbf{P}^{\mathrm{ex}};\ P^{\mathrm{cap}};\right]'$ with $\mathbf{B}$, $\mathbf{P}^{\mathrm{B}}$, $\mathbf{P}^{\mathrm{ex}}$ the vectors for $B_{i,s}(t), P_{i,s}^{\mathrm{B}}(t), P_{i,s}^{\mathrm{ex}}(t)$; $\mathbf{A}$, $\mathbf{C}$ and $\mathbf{D}$ are parameter matrixes.

According to the duality theory [39], the corresponding dual problem formulation is as follows:

$$\max_{\mathbf{Y}} \mathbf{D}' \cdot \mathbf{Y} \tag{36}$$

$$\text{s.t. } \mathbf{Y}' \cdot \mathbf{A} \le \mathbf{C}' \tag{37}$$

where $\mathbf{Y}$ is the vector of dual variables.

The dual problem of (Pb-1) can be explicitly written as the following (Db-1). The objective function (38) is at the bottom of this page.

$$-\lambda_{i,s}(t+1) + \lambda_{i,s}(t) - \gamma_{i,s}^{\mathrm{B}}(t) + v_{i,s}^{\mathrm{B}}(t) = 0 \quad \forall t \in \{1,2,...,T-1\} \tag{39}$$

$$\lambda_{i,s}(t) - \gamma_{i,s}^{\mathrm{B}}(t) + v_{i,s}^{\mathrm{B}}(t) = 0 \quad \text{when } t = T \tag{40}$$

$$-\lambda_{i,s}(t)\Delta t + \beta_{i,s}(t) - \gamma_{i,s}^{\mathrm{P}}(t) + v_{i,s}^{\mathrm{P}}(t) = 0 \tag{41}$$

$$-\beta_{i,s}(t) - \theta_s(t) + \eta_s(t) - \varepsilon_{i,s}(t) + \sigma_{i,s}(t) = \tau_s C(t)\Delta t \tag{42}$$

$$\sum_{s=1}^{S}\sum_{t=1}^{T}\theta_s(t) + \sum_{s=1}^{S}\sum_{t=1}^{T}\eta_s(t) = \rho \tag{43}$$

$$\{\gamma_{i,s}^{\mathrm{B}}(t), \gamma_{i,s}^{\mathrm{P}}(t), \theta_s(t), v_{i,s}^{\mathrm{B}}(t), v_{i,s}^{\mathrm{p}}(t), \eta_s(t), \varepsilon_{i,s}(t), \sigma_{i,s}(t)\} \ge 0 \tag{44}$$

$\gamma_{i,s}^{\mathrm{P}}(t)$ and $v_{i,s}^{\mathrm{p}}(t)$ are dual variables to the upper and lower bounds of aggregate battery charging/discharging power in constraint (9); $\varepsilon_s(t)$ and $\sigma_s(t)$ are dual variables to the upper and lower bounds of power flow in constraint (10). $\theta_s(t)$ and $\eta_s(t)$ are dual variables to the upper and lower bounds of power exchange in constraint (11).

$$\min_{\substack{\lambda_{i,s}(t);\gamma_{i,s}^{\mathrm{B}}(t);v_{i,s}^{\mathrm{B}}(t);\\ \gamma_{i,s}^{\mathrm{p}}(t);v_{i,s}^{\mathrm{p}}(t);\theta_s(t);\eta_s(t)}} \sum_{s=1}^{S}\sum_{i=1}^{I}\sum_{t=1}^{T} \begin{aligned}\{&B_{i,s}(0)\cdot\lambda_{i,s}(1) + \Delta_{i,s}^{B}(t)\cdot\lambda_{i,s}(t) + (P_{i,s}^{\mathrm{W}}(t) - P_{i,s}^{\mathrm{L}}(t))\beta_{i,s}(t) - \overline{B}_{i,s}(t)\cdot\gamma_{i,s}^{\mathrm{B}}(t)\\ &+ \underline{B}_{i,s}(t)\cdot v_{i,s}^{\mathrm{B}}(t) - \overline{P}_{i,s}(t)\cdot\gamma_{i,s}^{\mathrm{p}}(t) + \underline{P}_{i,s}(t)\cdot v_{i,s}^{\mathrm{p}}(t) - \overline{P}_i^{\mathrm{ex}}\cdot\varepsilon_{i,s}(t) - \overline{P}_i^{\mathrm{ex}}\cdot\sigma_{i,s}(t)\}\end{aligned} \tag{38}$$

According to the duality theory [39], the optimal variable is interpreted as shadow price, or the improvement in the objective function value versus the per unit relaxation of the associated constraint. We firstly consider a single scenario problem (sinPb-1):

$$\min_{P_{i,s}^{\text{ex}}(t),P^{\text{cap}}} \sum_{i=1}^{I}\sum_{t=1}^{T} C(t)\cdot P_{i,s}^{\text{ex}}(t)\Delta t+\rho\cdot P^{\text{cap}} \quad (45)$$
$$\text{s.t. (6)-(11)}$$

Let $\theta_s^*(t)$ and $\eta_s^*(t)$ be the shadow price of constraint (11), then the value of $\theta_s^*(t)$ (or $\eta_s^*(t)$) is the reduction in function value if the upper bound of total power exchange $P^{\text{cap}}$ (or lower bound of total power exchange $-P^{\text{cap}}$) is increased (or decreased) by a unit of 1. With the function $U(P_{i,s}^{\text{ex}}(t))=\sum_{i=1}^{I}\sum_{t=1}^{T}C_t\cdot P_{i,s}^{\text{ex}}(t)\Delta t$, the following equations determine $\theta_s^*(t)$ and $\eta_s^*(t)$ [39]:

$$\theta_s^*(t)=\left(\frac{\partial U(P_{i,s}^{\text{ex}*}(t))}{\partial(P^{\text{cap}})}\right)/\left(\frac{\partial(\sum_{i=1}^{I}P_{i,s}^{\text{ex}*}(t))}{\partial(P^{\text{cap}})}\right)=\frac{\Delta U(P_{i,s}^{\text{ex}*}(t))}{\Delta(\sum_{i=1}^{I}P_{i,s}^{\text{ex}*}(t))} \quad (46)$$

$$\eta_s^*(t)=\left(\frac{\partial U(P_{i,s}^{\text{ex}*}(t))}{\partial(P^{\text{cap}})}\right)/\left(\frac{\partial(-\sum_{i=1}^{I}P_{i,s}^{\text{ex}*}(t))}{\partial(P^{\text{cap}})}\right)=\frac{\Delta U(P_{i,s}^{\text{ex}*}(t))}{\Delta(-\sum_{i=1}^{I}P_{i,s}^{\text{ex}*}(t))} \quad (47)$$

where $P_{i,s}^{\text{ex}*}(t)$ comprises the optimal solution to (sinPb-1). Equation (46) and (47) illustrate a fact that the value of $\theta_s^*(t)-\eta_s^*(t)$ represents the relative change of the objective function due to the change in the constraint (11) given a relaxation of $P^{\text{cap}}$. Intuitively, a relatively large value of $\theta_s^*(t)-\eta_s^*(t)$ suggests a great value of increasing $P^{\text{cap}}$, which implies that at moment $t$ the demand for electricity is great. By strict complementarity [39], the value $\theta_s^*(t)-\eta_s^*(t)$ can illustrate at which moment the total power exchange has reached its upper bound $P^{\text{cap}}$ ($\theta_s^*(t)-\eta_s^*(t)>0$ since $\theta_s^*(t)>0,\ \eta_s^*(t)=0$) or its lower bound $-P^{\text{cap}}$ ($\theta_s^*(t)-\eta_s^*(t)<0$ since $\theta_s^*(t)=0,\ \eta_s^*(t)>0$). Therefore, the values of $\theta_s^*(t)-\eta_s^*(t)$ represent the value (or marginal price) of the total electricity exchange during $t$. Then we formulated the following problem (sinPb-2):

$$\min_{P_{i,s}^{\text{ex}}(t)} \sum_{i=1}^{I}\sum_{t=1}^{T} C_s^*(t)\cdot P_{i,s}^{\text{ex}}(t)\Delta t \quad (48)$$
$$\text{s.t. (6)-(10)}$$

where $C_s^*(t)=C_t+\theta_s^*(t)-\eta_s^*(t)$, $\theta_s^*(t)$ and $\eta_s^*(t)$ are the optimal values of the dual variables to the primal problem (sinPb-1).It is demonstrated that the solution of problem (sinPb-2) also solves problem (sinPb-1) by proving that optimal conditions and the feasibility conditions are equivalent in two problems formulations [40]. It should be pointed out that the problem (sinPb-2) can be further decomposed into $I$ subproblems, since the constraint (10) which couples multiple microgrids is relaxed in problem (sinPb-2) and each microgrid can carry out its own management. Following this line, we try to handle the stochastic problem (Pb-1). We can solve (Db-1) and obtain $S$ sets of dual variable vectors, considering uncertainties of EVs and wind power. Among these optimal values of dual variables, $\theta_s^*(t)$ and $\eta_s^*(t)$ are the shadow price of constraint (11). Then we utilize these two variables to generate the updated price signal. Since the probability of each scenario ($\tau_s$) is considered in the objective function of stochastic formulation (Pb-1), $\theta_s^*(t)$ and $\eta_s^*(t)$ contain the probability information. At each moment $t$, $\theta_s^*(t)-\eta_s^*(t)\neq 0$ means that it reaches the bound with the probability $\tau_s$.

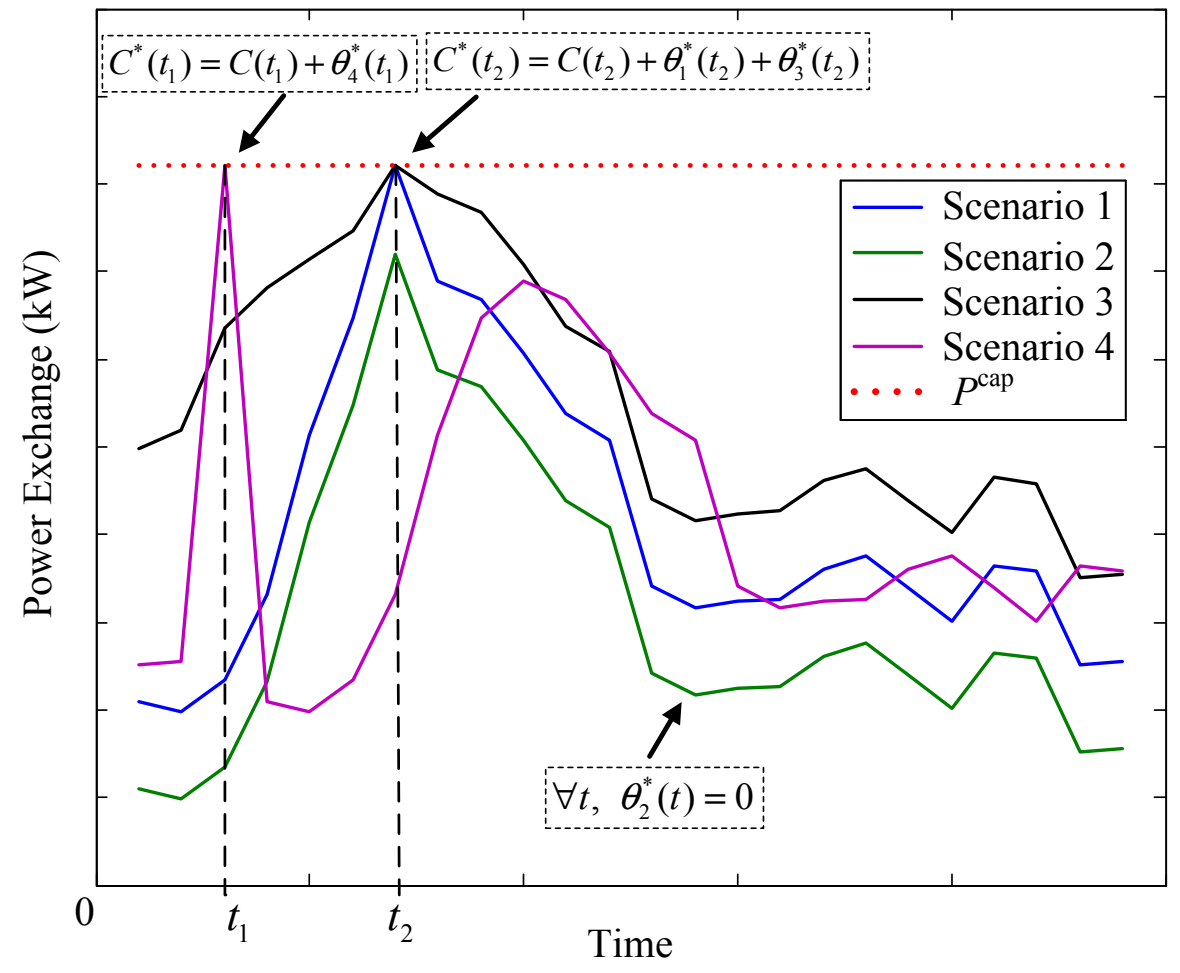

Fig.2. Updated price signal for different scenarios

A four scenario example is shown in Fig. 2 to give a further explanation about the pricing scheme. Since the power exchange in these four scenarios cannot reach the lower bound, $\forall t,s\ \eta_s^*(t)=0$ and we only need to consider the dual variable $\theta_s^*(t)$ which is related to the upper bound. The power exchange in scenario 2 cannot reach the upper bound. Thus, $\forall t,\ \theta_2^*(t)=0$. Scenario 1, 2 and 3 can all be seen as worst cases, because at some moments the power exchange reaches the limits. For scenario 1 and scenario 3, they reach the same upper bound of the power exchange at $t_2$. We can obtain that $\begin{cases}\theta_1^*(t_2)\neq 0,\theta_3^*(t_2)\neq 0\\ \theta_1^*(t)=\theta_3^*(t)=0\quad \forall t\neq t_2\end{cases}$. For scenario 4, the power exchange reaches the limit at $t_1$. We can also get $\theta_4^*(t_1)\neq 0$. According to (43), we can conclude that

$$\begin{aligned}&\sum_{s=1}^{S}\sum_{t=1}^{T}\theta_s^*(t)+\sum_{s=1}^{S}\sum_{t=1}^{T}\eta_s^*(t)=\theta_1^*(t_2)+\theta_3^*(t_2)+\theta_4^*(t_1)=\rho\\&C^*(t_1)=C(t_1)+\theta_4^*(t_1)\\&C^*(t_2)=C(t_2)+\theta_1^*(t_2)+\theta_3^*(t_2)\end{aligned} \quad (49)$$

For the system studied in this paper, we assume that wind power is not able to reach a very high penetration level that will drastically change the power exchange profiles. So there are only one or a few scenarios in which the exchange power can reach the limits. We can guarantee that the power exchange within the safe range when $\sum_{s=1}^{S}\theta_s^*(t)-\sum_{s=1}^{S}\eta_s^*(t)$ is utilized as an adjustment to the original TOU price. With this updated price signal $C^*(t) = C(t) + \sum_{s=1}^{S}\theta_s^*(t) - \sum_{s=1}^{S}\eta_s^*(t)$, we can get the decentralized formulation (Db-2).

$$\min_{P_{i,s}^{\mathrm{ex}}(t)} W = \sum_{s=1}^{S}\sum_{i=1}^{I}\sum_{t=1}^{T}\tau_s C^*(t)\cdot P_{i,s}^{\mathrm{ex}}(t)\Delta t \tag{50}$$
s.t. (6)-(10)

where $W$ is the average electricity cost of (Db-2). To make this value equal to the MGCC payment for exchanged energy to the main grid, we introduce a scaling factor $\varepsilon = W'/W$.

$$W' = \sum_{s=1}^{S}\sum_{i=1}^{I}\sum_{t=1}^{T}\tau_s C(t)\cdot P_{i,s}^{\mathrm{ex}*}(t)\Delta t \tag{51}$$

where $P_{i,s}^{\mathrm{ex}*}(t)$ is the optimal values of $P_{i,s}^{\mathrm{ex}}(t)$ in (Db-2). Then we can get the final updated price signal:

$$\overline{C^*}(t) = [C^*(t) + \sum_{s=1}^{S}\theta_s^*(t) - \sum_{s=1}^{S}\eta_s^*(t)]\cdot\varepsilon \tag{52}$$

It will finally be broadcasted to each microgrid and the decentralized scheduling of multi-microgrid system is described as follows:

**Algorithm 1**:***Decentralized Scheduling Strategy***

1. Obtain each microgrid's information on load profile $P_i^{\mathrm{L}}(t)$, probability $\tau_s$, wind power $P_{i,s}^{\mathrm{W}}(t)$ and battery storage $\underline{P}_{i,s}(t)$ $\overline{P}_{i,s}(t)$ $\underline{B}_{i,s}(t)$ $\overline{B}_{i,s}(t)$ //using the scenario generation and reduction in Section IV.
2. Solve the dual problem (Db-1) and obtain the optimal dual variables $\theta_s^*(t)$ $\eta_s^*(t)$.
3. Solve the problem (Db-2) and obtain the optimal value of $P_{i,s}^{\mathrm{ex}*}(t)$ and objective function $W$.
4. Calculate the MGCC payment $W'$ for the energy exchange with the main grid.
5. Calculate the updated price signal as:
$$\overline{C^*}(t) = [C^*(t) + \sum_{s=1}^{S}\theta_s^*(t) - \sum_{s=1}^{S}\eta_s^*(t)]\cdot(W'/W)$$
and broadcast the control signal to all microgrids.
6. **For** Microgrid $i$←1 to $I$ **Do**
   optimize its own battery charging/discharging power to minimize the electricity cost $\sum_{t=1}^{T}\overline{C^*}(t)\cdot P_i^{\mathrm{ex}}(t)\Delta t$, subjected to the constraints (6)-(10).
7. **End for**

Fig.3 shows the information exchange between the MGCC and each microgrid for the proposed algorithm. Given the updated price signal broadcasted by the MGCC, each microgrid adjusts its own charging/discharging profiles independently. Although the MGCC need to get the global information, it will not carry out specific control or scheduling for each individual microgrid operation. It should be pointed out that the only constraint (11) which couples multiple microgrids is relaxed by adopting the updated price signal. So each microgrid can carry out its own strategy. It just needs to solve the problem as the Step 6 of Algorithm 1, without considering the interaction of other microgrids.

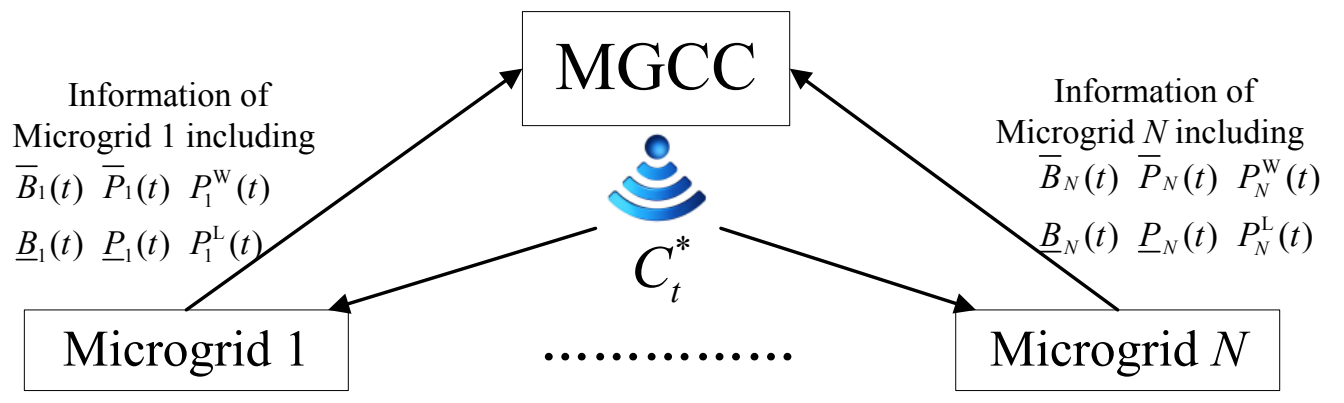

Fig. 3. Schematic view of the information flow between the MGCC and microgrids. Given the updated price signal, microgrids update their operations independently. The MGCC guides their updates by altering the price signal.

## VI. Numerical Testing Results

Numerical test of the proposed method is performed in a five microgrids system, as shown in Fig. 4. All microgrids are connected to a 400V bus. The schedule period starts from 12:00 to 12:00 next day. The baseline load profile of each microgrid covering the power consumption of user general demand without EVs over 24 hours is shown in Table I.

TABLE I
Load profiles of each microgrid

| Hour | Load(kW) | | | | |
|---|---|---|---|---|---|
| | MG1 | MG2 | MG3 | MG4 | MG5 |
| 13 | 953 | 636 | 95 | 185 | 216 |
| 14 | 953 | 636 | 225 | 420 | 516 |
| 15 | 947 | 631 | 215 | 430 | 516 |
| 16 | 967 | 644 | 352 | 743 | 876 |
| 17 | 953 | 636 | 563 | 1132 | 1356 |
| 18 | 1053 | 702 | 655 | 1340 | 1596 |
| 19 | 1033 | 689 | 1369 | 2681 | 2160 |
| 20 | 1000 | 667 | 890 | 1780 | 1643 |
| 21 | 967 | 644 | 865 | 1760 | 1400 |
| 22 | 820 | 547 | 827 | 1648 | 1320 |
| 23 | 700 | 467 | 624 | 1251 | 1500 |
| 0 | 593 | 396 | 266 | 559 | 660 |
| 1 | 367 | 244 | 86 | 139 | 180 |
| 2 | 353 | 236 | 47 | 103 | 120 |
| 3 | 333 | 222 | 79 | 144 | 180 |
| 4 | 333 | 231 | 68 | 127 | 147 |
| 5 | 433 | 289 | 74 | 151 | 180 |
| 6 | 387 | 261 | 72 | 150 | 180 |
| 7 | 520 | 347 | 38 | 67 | 84 |
| 8 | 567 | 378 | 98 | 202 | 240 |
| 9 | 820 | 547 | 114 | 216 | 264 |
| 10 | 1053 | 702 | 90 | 151 | 191 |
| 11 | 967 | 654 | 77 | 147 | 170 |
| 12 | 973 | 644 | 75 | 150 | 258 |

Each microgrid has a wind turbine generator, whose rated power is set as 500kW. Other parameters, $v_{\mathrm{ci}}$, $v_r$, $v_{\mathrm{co}}$ are set as 3m/s, 10m/s, 20m/s respectively. The data on forecasted wind speed are collected from Tianjin Eco-city [41]. It is assumed that all wind power generators produce active power at a unity power factor, i.e., neither requesting nor producing reactive power. The TOU price has three stages, as shown in Table II.

In the simulation, the total number of EVs in each microgrid is 100, and all EVs are with the same specification where the battery capacity is 33kWh to ensure that the vehicle can finish daily driving mission. The initial energy of EVs is uniformly distributed between 10% and 50% of the battery capacity. To avoid excessive discharge and damage to the battery, the SOC is required to be always beyond 10%. Electric drive efficiency is 6.7 km/kWh [42].

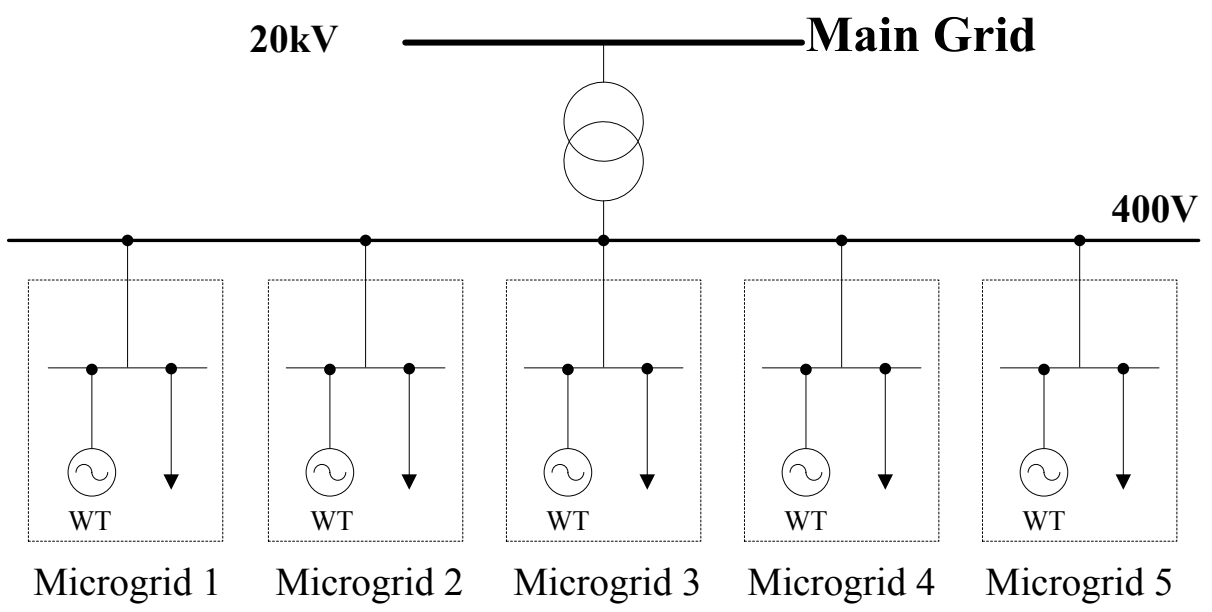

Fig.4. Study case multi-microgrid system.

To represent the stochastic nature of the wind power and EVs, Monte Carlo method is used to generate a total number of 2000 scenarios according to the statistical distributions described in Section IV. The fast forward algorithm [33] is then utilized to reduce the number of scenarios to 100 in light of the computational burden.

TABLE II
TOU Electricity Price

| Period | Price(CNY/kWh) |
|---|---|
| 0:00-7:00 | 0.3515 |
| 7:00-12:00 | 0.8135 |
| 12:00-19:00 | 0.4883 |
| 19:00-23:00 | 0.8135 |
| 23:00-24:00 | 0.3515 |

## A. *Parameters of the Aggregate Battery Storage*

The aggregate battery storage consists of all EVs connected to the grid. Due to the individual behaviors of each EV, the parameters of the aggregate battery storage, including $\underline{B}_{i,s}(t)$, $\overline{B}_{i,s}(t)$, $\underline{P}_{i,s}(t)$ and $\overline{P}_{i,s}(t)$, cannot be considered as constant values. We firstly investigate the influence of driving patterns on the performance of the aggregate battery storage. To better illustrate this influence, we select the scenario with the highest probability to show the variability of the aggregate battery storage. The bound of charging/discharging power and the storage capacity in this scenario are shown in Table III. It can be seen that these two parameters fluctuate during the schedule period. They are determined by the arrival time, departure time and driving distance. The maximum charging/discharging power and the upper bound of battery energy are in proportion to the number of EVs at home. The lower bound of battery remaining energy is determined by the number of EVs plugged in and their driving distance. Each EV must leave home with enough energy for its driving missions. From 20:00 to 6:00, since most of EVs are at home and can be plugged into the grid, there is a high storage capacity and wide range of charging/discharging power. Conversely, from 10:00 to 15:00, most of EVs are not at home. In this case, the aggregate battery storage can rarely contribute to the energy scheduling, due to the relatively low storage capacity and charging/discharging power.

TABLE.III
Battery Information

| Time(h) | $n_t$ | $\overline{P}_{i,t}$ (kW) | $\overline{B}_{i,t}$ (kWh) | $\underline{B}_{i,t}$ (kWh) |
|---|---|---|---|---|
| 13 | 1 | 3 | 33 | 3 |
| 14 | 2 | 6 | 66 | 7 |
| 15 | 4 | 12 | 132 | 13 |
| 16 | 19 | 57 | 627 | 63 |
| 17 | 33 | 99 | 1089 | 109 |
| 18 | 45 | 135 | 1485 | 149 |
| 19 | 55 | 165 | 1815 | 182 |
| 20 | 67 | 201 | 2211 | 221 |
| 21 | 77 | 231 | 2541 | 254 |
| 22 | 89 | 267 | 2937 | 294 |
| 23 | 97 | 291 | 3201 | 320 |
| 24 | 100 | 300 | 3300 | 330 |
| 1 | 100 | 300 | 3300 | 330 |
| 2 | 100 | 300 | 3300 | 330 |
| 3 | 100 | 300 | 3300 | 330 |
| 4 | 100 | 300 | 3300 | 365 |
| 5 | 96 | 288 | 3168 | 371 |
| 6 | 86 | 258 | 2838 | 465 |
| 7 | 58 | 174 | 1914 | 475 |
| 8 | 21 | 63 | 693 | 169 |
| 9 | 4 | 12 | 132 | 31 |
| 10 | 0 | 0 | 0 | 0 |
| 11 | 0 | 0 | 0 | 0 |
| 12 | 0 | 0 | 0 | 0 |

## B. *Performance of the Proposed Coordination Strategy*

In the proposed strategy, coordinated charging/discharging of EVs is integrated into the hierarchical energy management framework. Working as energy storage devices, EVs help to lower the power exchange peak as well as save the total electricity cost. The profiles of the average power exchange in two cases are shown in Fig. 5. Uncoordinated strategy means that each EV just gets its required energy for driving, while not serving as mobile storage device. Since we take the uncertainty of wind power generation into consideration, the value of $P^{cap}$ is determined by the worst case. It can be seen that with coordinated charging/discharging scheduling, the power exchange peak is decreased by 10.18%, compared with the uncoordinated scheduling strategy.

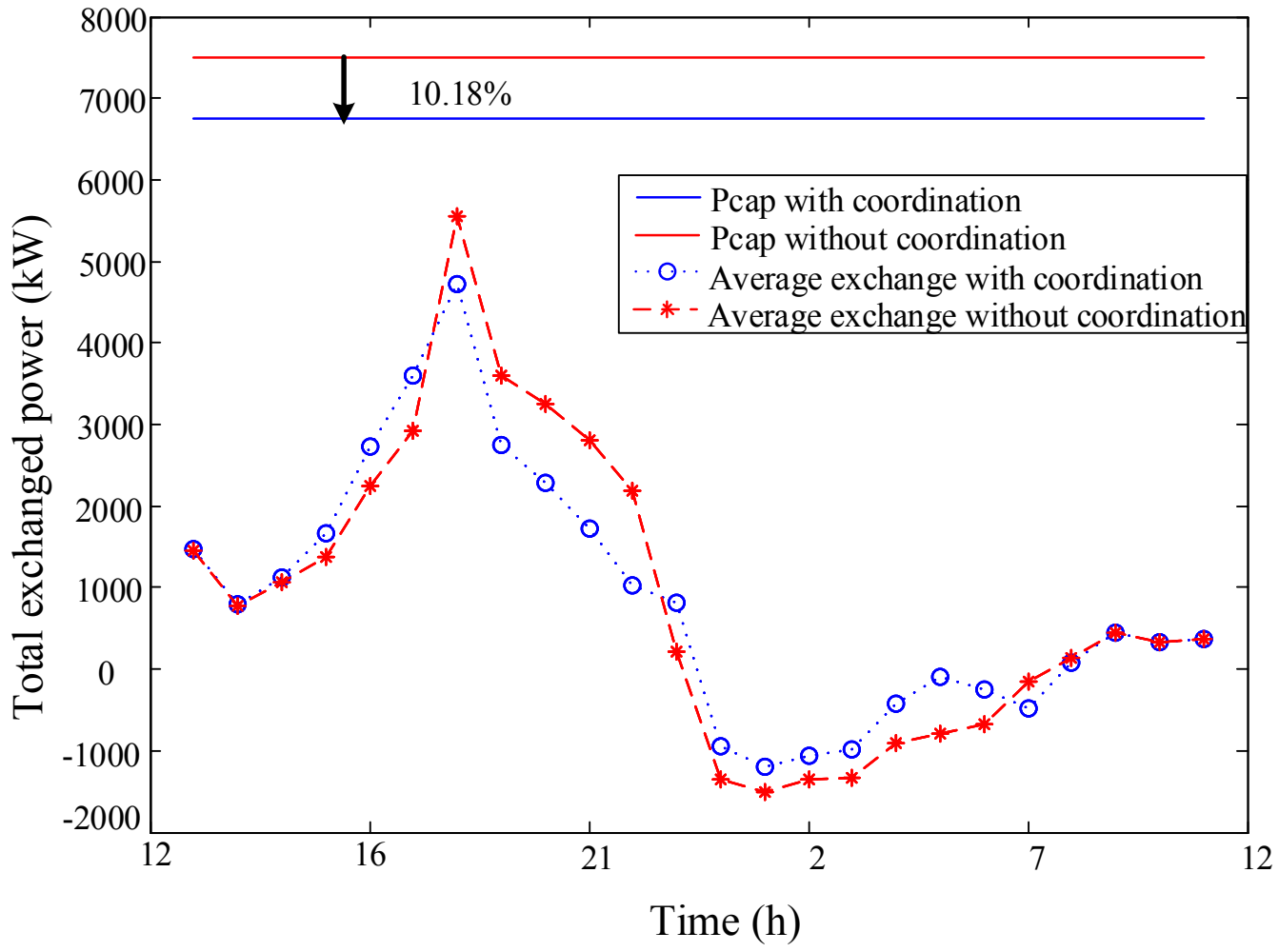

Fig. 5 Total power exchange with coordination and without coordination.

Besides, although the total amounts of energy exchange in both cases are the same, the average power exchange with coordination performs less volatility than that without coordination. We can see that the power exchange with coordination is lower during the price-peak time (19:00-23:00) and higher during the price-valley time (23:00-7:00) than that without coordination. To be explicit, the coordinated strategy makes the aggregate battery storage discharge during high price period and charge during low price period. Therefore, the total electricity cost will be decreased with the proposed method. Some detailed results of two cases are shown in Table IV. The value of $P^{cap}$ decreases from 7512kW to 6747. On average, the total electricity cost is reduced by 12.97% and the standard deviation (SD) is reduced by 16.49%.

TABLE.IV
Comparison of the power exchange profiles with coordination and without coordination

| | $P^{cap}$ (kW) | Average cost (CNY) | SD |
|---|---|---|---|
| Without coordination | 7512 | 15215 | 1850 |
| With coordination | 6747 | 13242 | 1545 |

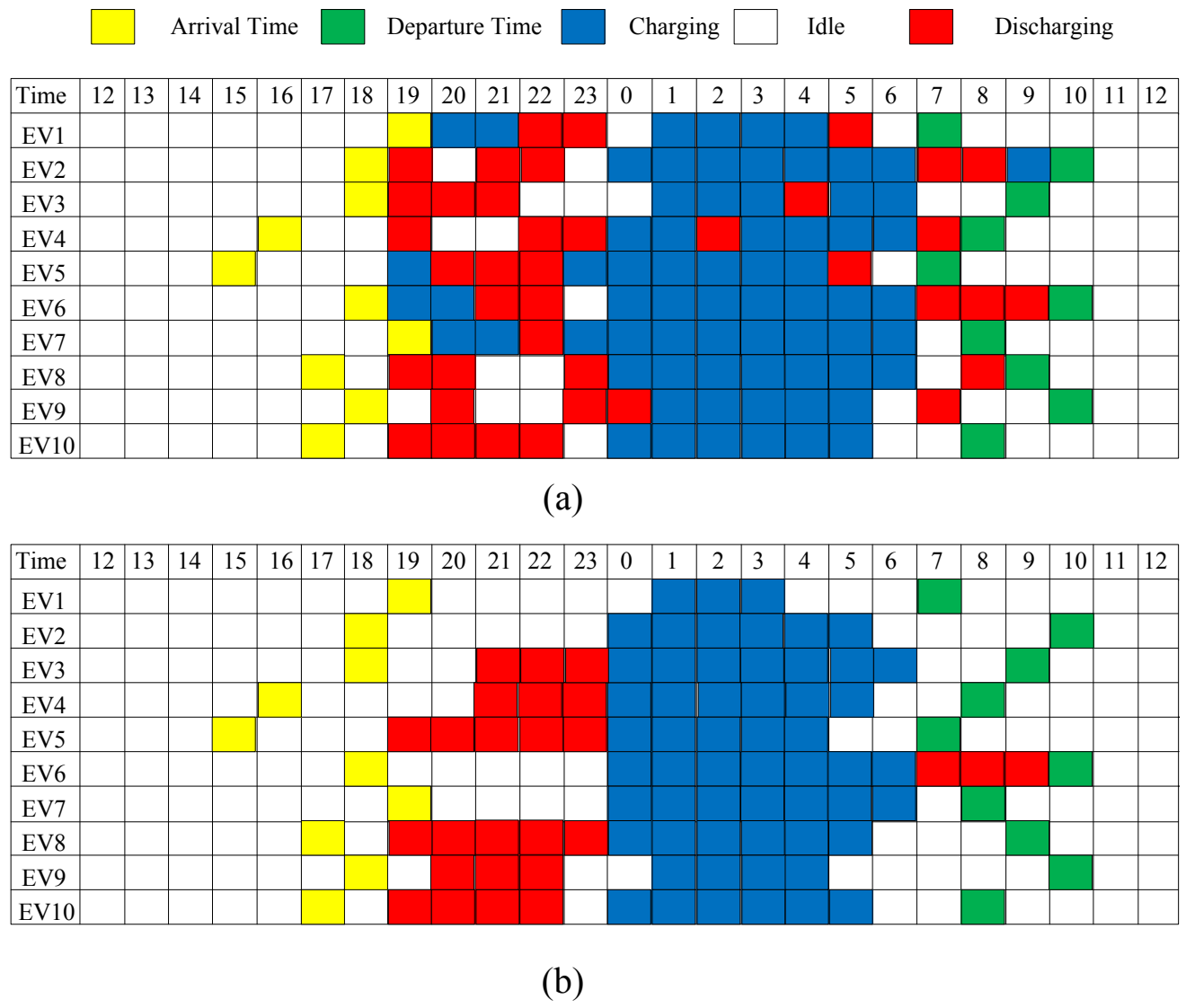

Fig. 6. Comparison of schedule results: (a) the MCLF allocation; (b) minimum state transitions allocation.

For the EV aggregator scheduling problem, the objective is to allocate the total charging/discharging power with the minimum cost. To give an intuitional demonstration of the allocation result, the scheduling results of EV 1-10 in the highest probability scenario are extracted and illustrated in Fig.6. We compare minimum state transitions allocation with the minimum-cost-least-fluctuation (MCLF) allocation [29]. MCLF allocation means that the aggregator schedules the charging/discharging of each EV, aiming to achieve the least fluctuation charging/discharging profiles under the premise the minimum cost. Minimum state transitions allocation refers to a strategy which is described in Section IV-C. As shown in Fig. 6(a), the charging and discharging states frequently change, which have great negative effects on batteries' cycle life. In Fig. 6(b), the situation has a great improvement. The proposed allocation method does not generate charging/discharging cycles more frequently than necessary, even if the battery cost is not considered in the MGCC scheduling stage. From the manual of "Nissan Leaf" [42], the cycle life of battery is 3000 and the battery price is 78,000 CNY. Under the simple assumption that each charging and discharging cycle causes the same wear on the battery, the penalty would be 26 CNY per cycle. For 100 EVs in a microgrid, the proposed allocation strategy saves 67 cycles on average. If the EV aggregator pays the penalty of battery cycle life to EV owners, the proposed method will save 1742 CNY for the aggregator. Moreover, if it is assumed that at least one cycle is required for driving demand, the proposed strategy does not introduce any extra cycle on batteries. Intuitively, the number of EVs connected to the grid is one of the important factors which will determine the scheduling results on battery degradation cost. When more EVs are connected to the grid, the aggregator has a more flexibility to choose the proper EVs to work as battery storages without introducing extra cycle cost. In this case study, the result shows that no extra cost is introduced until the number of EVs decreases from 100 to 89 in each microgrid.

### C. *Performance of the Decentralized Scheduling Strategy*

From problem (Pb-1), we can obtain the $P^{cap}$=6747kW, which is the upper bound of total power exchange between the multi-microgrid system and the main distribution grid. This ensures that the peak of power exchange would not have a great impact on the distribution grid. The average results of power exchange under the centralized strategy and the price-based decentralized strategy are compared in Fig. 7. Both of them have the same average total energy exchange during the scheduling horizon. For the centralized strategy, the MGCC perform the management according to (5)-(11) to achieve the minimum cost. For the price-based decentralized scheduling strategy, a peak price is added to the original TOU price, as shown in Fig. 7. The electricity price during 18:00-19:00 reaches 1.0580 CNY, which is much higher than the original 0.4883CNY. The updated price has four price stages, as in Table V. Responding to the updated price, each microgrid adjusts its charging/discharging schedule independently, instead of given by a centralized mechanism.

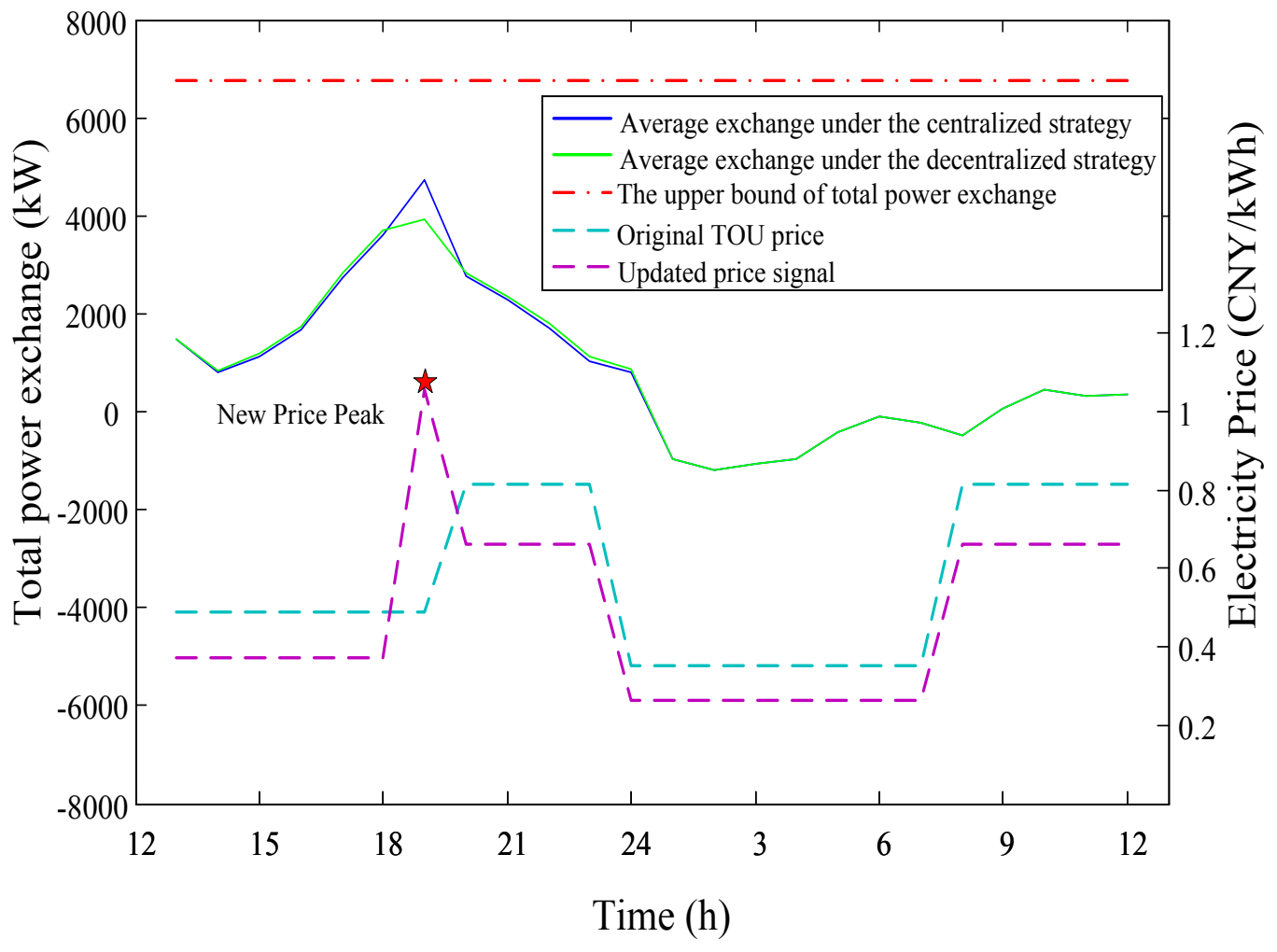

Fig. 7. Total power exchange under centralized and decentralized strategies.

TABLE V
Updated Price Signal

| Period | Price(CNY/kWh) |
|---|---|
| 0:00-7:00 | 0.2622 |
| 7:00-12:00 | 0.6619 |
| 12:00-18:00 | 0.3733 |
| 18:00-19:00 | 1.0580 |
| 19:00-23:00 | 0.6619 |
| 23:00-24:00 | 0.2622 |

It can be seen in Fig.7 that these two profiles are not totally overlapped. The average power exchange under the price-based decentralized strategy is higher during 12:00-18:00 and 19:00-23:00, and lower during 18:00-19:00 than that under the centralized strategy. It is because that to ensure the robustness, the price-based strategy adjusts the price at which moment there may be power exchange peak. It results in that microgrids in low power exchange scenarios also update their schedule to discharging the battery at 19:00 to save costs. Therefore, such amount of power exchange shifts to 12:00-18:00 and 19:00-23:00, and the profile is further leveled. It should be noted that the electricity cost of the proposed decentralized strategy is not cost-equivalent to that of the centralized strategy, due to the power shift described above. From Table VI, we see that the average cost of the proposed decentralized strategy is 13452CNY, which is a little more than the cost of the centralized one.

TABLE VI
Optimization results under different strategies

| Strategy | On average | | Worst Case | |
|---|---|---|---|---|
| | Cost(CNY) | Peak(kW) | Cost(CNY) | Peak(kW) |
| Centralized strategy | 13362 | 4722 | 32956 | 6747 |
| Decentralized strategy under updated price signal | 13452 | 3942 | 32956 | 6747 |
| Decentralized strategy under TOU | | | 32634 | 8306 |

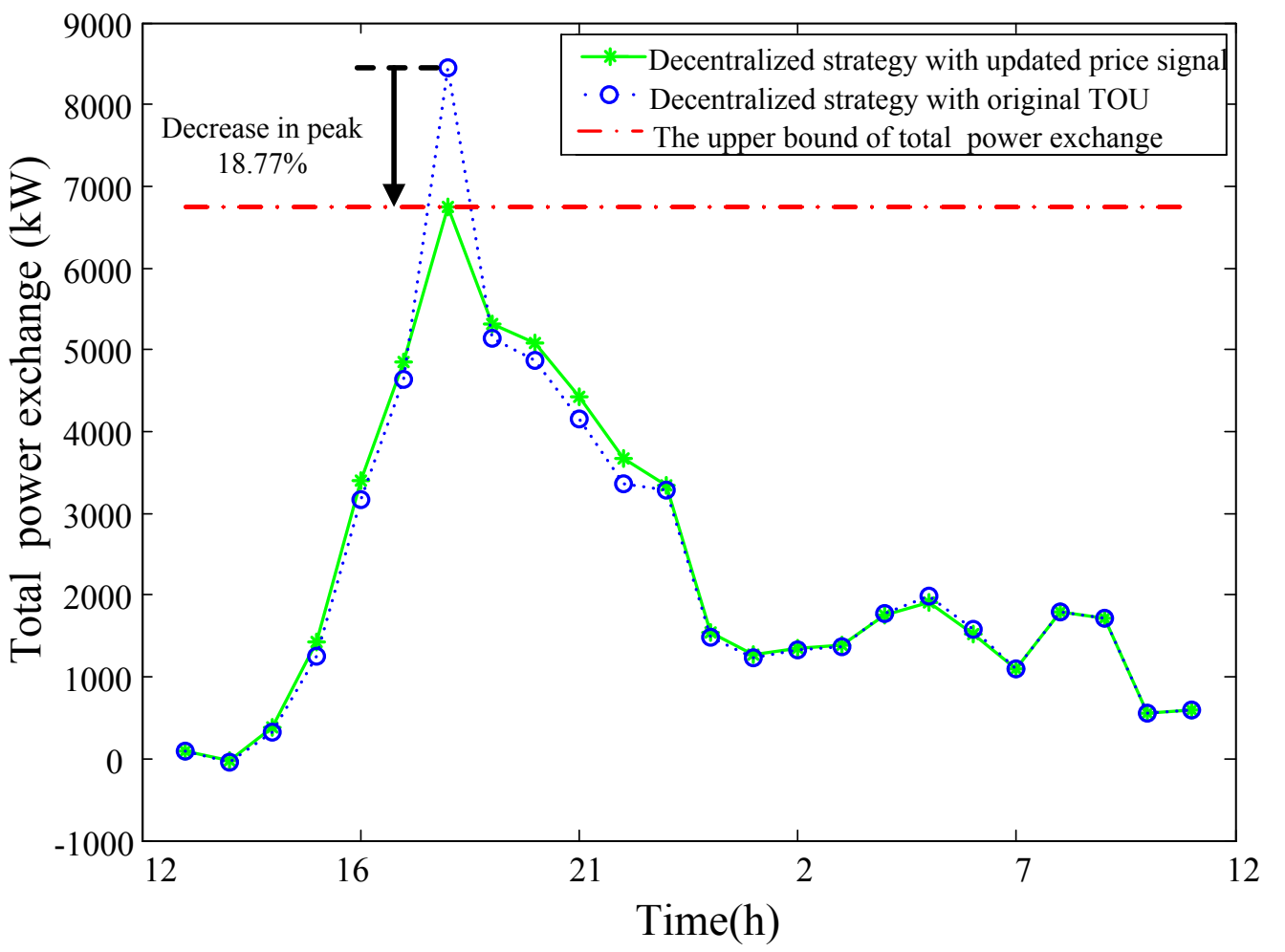

Fig. 8. Total power exchange under different strategies in the worst case.

Then, we investigate the worst case which has a great influence on the upper bound of power exchange. A comparison is made between the proposed price-based strategy and the decentralized strategy with the original TOU price, as shown in Fig. 8. The decentralized strategy with TOU price means that the constraint (11) is not considered and each microgrid makes decisions to minimize its own cost under TOU price regardless of others microgrids' behavior. With this strategy, although the multi-microgrid system can achieve the lowest cost, it cannot guarantee the power exchange within the safe range. We can see from Fig. 8 that the power exchange peak of the decentralized strategy under original TOU price is 8306kW, which will be hazardous for the power facilities and distribution network. As shown in Table VI, the peak of total exchange power drops to 6747kW, decreasing by 18.77% and staying within the $P^{cap}$. The centralized scheduling strategy proposed in Section IV obtains the same results as the proposed decentralized strategy. Since the MGCC has to limit the peak of power exchange lower than $P^{cap}$ to ensure the security of distribution network, the operation costs of both centralized and price-based decentralized strategies are increased by 0.98%, compared with the cost of decentralized strategy under TOU price.

To clearly show how each microgrid adjusts the aggregate battery charging/discharging power based on the updated price signal, we further investigate the performance of the aggregate battery storage in the worst case. From Fig. 9, it can be seen that under the updated price, the microgrid adjusts the battery charging/discharging schedule. From 12:00 to 18:00, the charging power is increased to the upper bound to get more energy. For individual scheduling under original TOU price, the total charging power is 146.8kW during 18:00-19:00, while for the updated price signal, the battery shifts to discharging at its maximum power (-165kW). The reason lies in the fact that a high peak price is at this moment and a substantial payback can be obtained if the battery is discharged. Therefore, the microgrid utilizes all EVs as microsources to feed power back to the grid during this period. From 19:00-23:00, the price is relatively high for both original TOU price and updated price signal. The battery storage would be discharged as much as possible. However, since some energy has been fed back during 18:00-19:00, the discharging power under updated price signal is lower than it under the original TOU from 19:00 to 23:00, as shown in Fig.9. After the price-peak time, the battery has the same amount of energy for both two strategies. In the remaining time, the charging/discharging schedules keep the same.

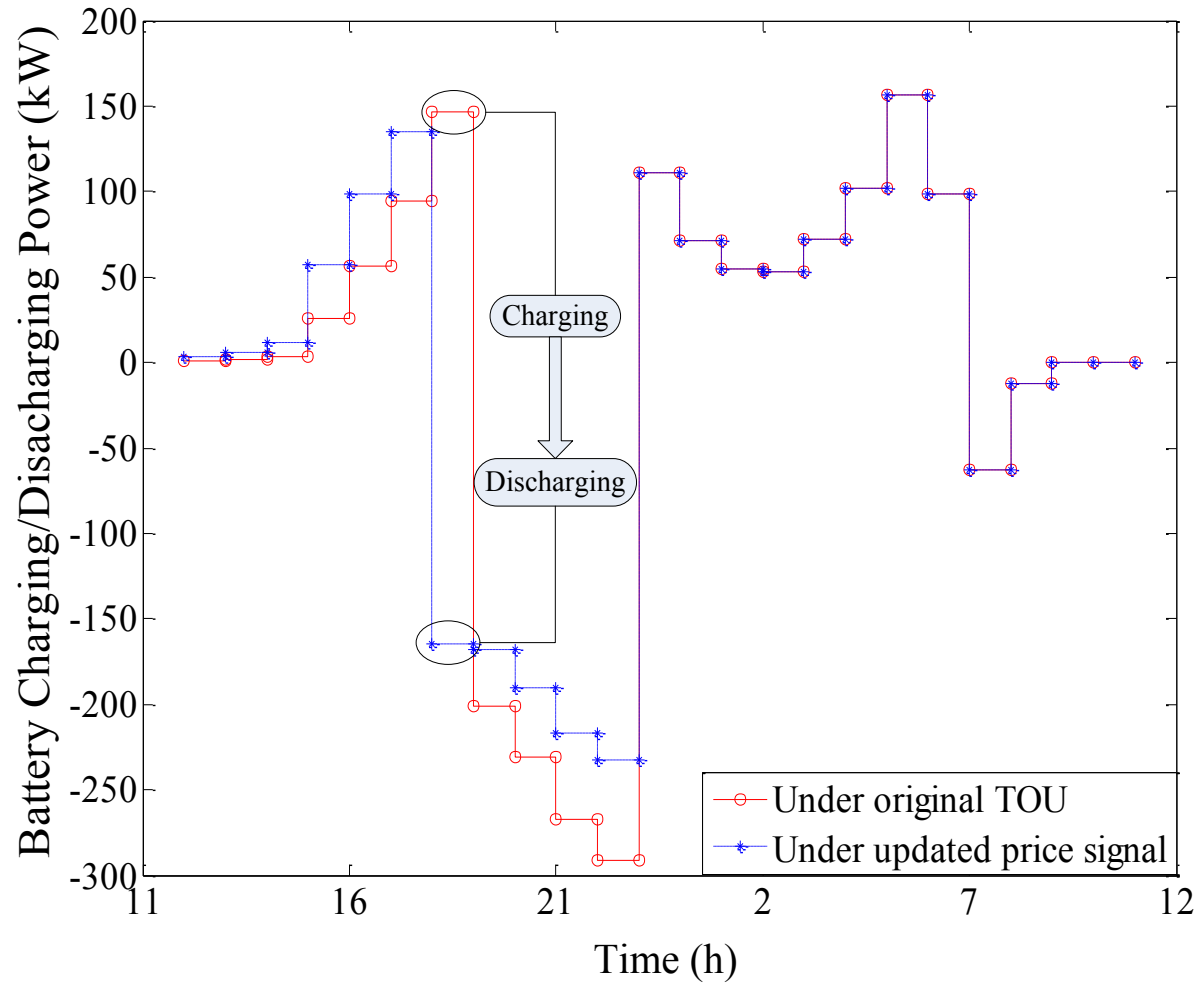

Fig. 9. Battery charging/discharging schedule under the original TOU price and the updated price signal in the worst case.

## VII. Conclusions and Future Work

This paper presents a coordinated scheduling framework and a two-stage optimization model for multi-microgrid system, considering EVs as energy storage devices. The first stage is MGCC scheduling stage, which tries to achieve the minimum cost by scheduling the total charging/discharging power of the aggregate battery storage under TOU price. The second stage is EV aggregator scheduling stage, which aims at allocating the total charging/discharging power obtaining from the first stage to each individual EV, and minimizing the transitions of battery charging/discharging states to decrease the wear on batteries. Scenario tree method is applied to solve the stochastic optimization problem with the uncertainties of wind power and EV charging/discharging behaviors. From the dual problem of the linear programming model, a price-based decentralized scheduling strategy for MGCC is presented based on the dynamic price update. The conclusions are drawn as follows:

1) The two-stage model proposed in this paper can effectively schedule the operation of the multi-microgrid system and utilize EVs as storage devices. The results show that the coordinated method restrains the power exchange peak within the safe range and decreases the wear on batteries caused by frequent transitions between charging/discharging states. In addition, the electricity cost is reduced by 12.97%, compared with the scheduling without EV coordination.
2) This paper presents a decentralized scheduling strategy based on dual variables' interpretation as shadow price. Simulation results show that the price-based strategy can effectively limit the total power exchange within the safe range. The power exchange peak is decreased by 18.77%, compared with the decentralized strategy under the original TOU price.

In the future work, we will incorporate the detailed AC operational constraints and system dynamics into consideration.

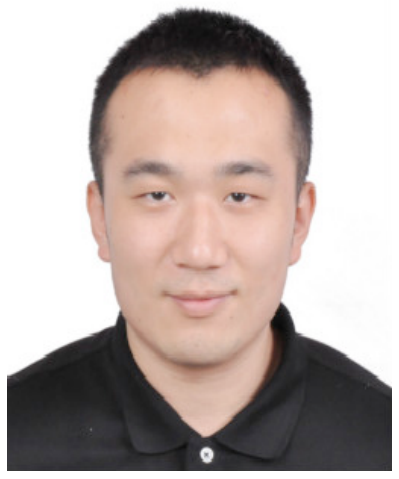

**Dai Wang** (S'12) received his B.S. degree in Electrical Engineering from Xi'an Jiaotong University, Shaanxi, China in 2010. He is currently a Ph.D. student at Systems Engineering Institute, Xian Jiaotong University and a visiting student at the Department of Mechanical Engineering, University of California, Berkeley. His research interests include demand response, electric vehicle, energy storage and electricity market.

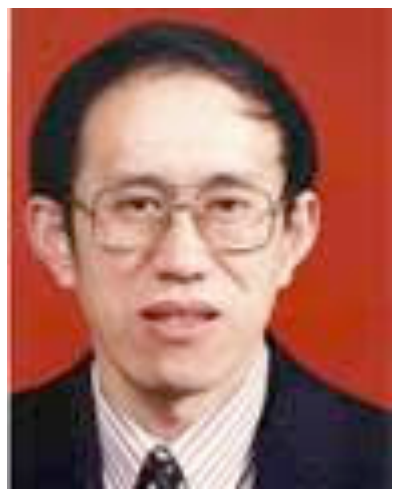

**Xiaohong Guan** (M'93-SM'95-F'07) received the B.S. and M.S. degrees in automatic control from Tsinghua University, Beijing, China, in 1982 and 1985, respectively, and the Ph.D. degree in electrical engineering from the University of Connecticut, Storrs, in 1993. He was a Consulting Engineer with PG&E from 1993 to 1995. From 1985 to 1988, and since 1995, he has been with the Systems Engineering Institute at Xi'an Jiaotong University, Xian, China, and currently he is the Cheung Kong Professor of systems engineering and Director of Systems Engineering Institute. He is also the Director of Center for Intelligent and Networked Systems, Tsinghua University. He visited the Division of Engineering and Applied Science, Harvard University, Cambridge, MA, from January 1999 to February 2000. His research interests include scheduling of power and manufacturing systems, bidding strategies for deregulated electric power markets, and economy and security of complex network systems.

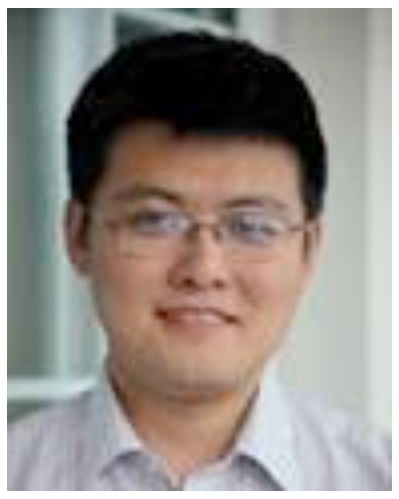

**Jiang Wu** (S'07-M'09) received his B.S. degree in Electrical Engineering from Xi'an Jiaotong University, Shaanxi, China in 2002, and his Ph.D. in Systems Engineering from the same university in 2008. He is currently an assistant professor at Systems Engineering Institute, Xian Jiaotong University. His research interests include stochastic unit commitment, deregulated electricity markets, and smart grid.

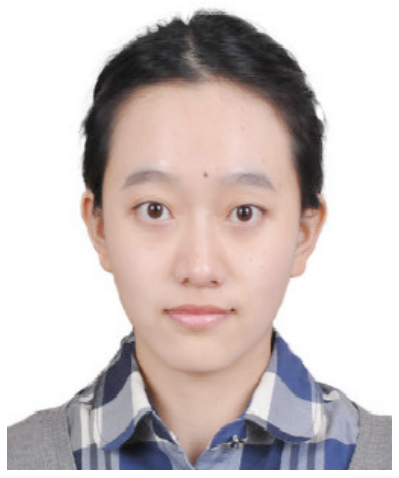

**Pan Li** received her B.S. degree in Electrical Engineering from Xi'an Jiaotong University, Shaanxi, China in 2011, and her M.S. degree in Systems Engineering from the same university in 2014. She is currently a Ph.D. student at the Clean Energy Institute, University of Washington, Seattle.
Her research interests include microgrid, solar power integration and smart grid.

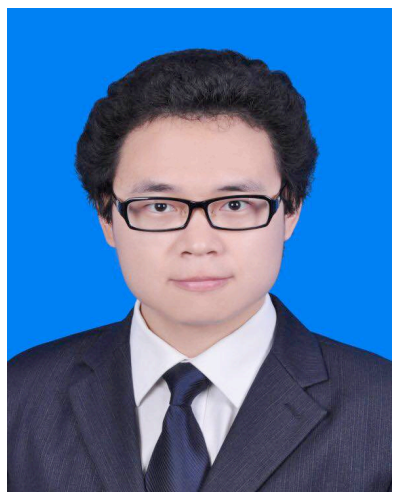

**Peng Zan** received his B.S. degree in Electrical Engineering from Xi'an Jiaotong University, Shaanxi, China in 2014. He is currently a Ph.D. student in the University of Maryland, College Park, MD. His research interests include demand response and smart grid.

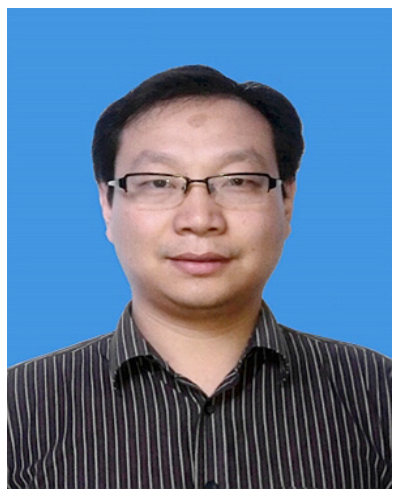

**Hui Xu** received his B.S. degree in Electrical Engineering from Xi'an Jiaotong University, Shaanxi, China in 2004, and his M.S. degree in Systems Engineering from the same university in 2008. He is currently with Binhai Power Supply Bureau, Tianjin, China. His research interests include distributed generation, distribution system automation, and power system planning.